\documentclass[12pt,reqno]{amsart}
\usepackage{fullpage}
\usepackage{amsthm,amsmath,amsfonts,amssymb,euscript,hyperref,graphics,color,slashed,mathrsfs,tikz,wrapfig,float,caption,enumerate,graphicx,mathtools,extarrows,cite,framed,mfirstuc}
\allowdisplaybreaks

\hypersetup{colorlinks, citecolor=blue, filecolor=blue, linkcolor=blue, urlcolor=blue}
\newcommand{\varp}{\varepsilon}
\newcommand{\pt}{\partial_t}
\newcommand{\px}{\partial_x}
\newcommand{\py}{\partial_y}
\newcommand{\pz}{\partial_z}
\newcommand{\Ga}{\Gamma}
\newcommand{\bu}{\bold{u}}
\newcommand{\norm}[1]{\|#1\|}
\newcommand{\les}{\lesssim}

\newcommand{\hcom}{H_{co}^m}
\newcommand{\hcoml}{H_{co}^{m-1}}

\newcommand{\zb}{Z^\beta}
\newcommand{\into}{\int_{\Omega}}
\newcommand{\curl}{\textup{curl}}
\numberwithin{equation}{section}

\newtheorem{Theorem}{Theorem}[section]
\newtheorem{Lemma}[Theorem]{Lemma}
\newtheorem{Proposition}[Theorem]{Proposition}

\newtheorem{Remark}[Theorem]{Remark}

\title[Zero-viscosity Limit for Anisotropic Boussinesq Equations]{Zero-viscosity Limit for  Boussinesq Equations with Vertical Viscosity and Navier Boundary in the Half Plane }
\author{Mengni Li, Yan-Lin Wang\textsuperscript{*}} 
\thanks{* Corresponding author}

\address{Mengni Li, School of Mathematics, Southeast University, Nanjing 211189, China}
\email{krisymengni@163.com}
\address{Yan-Lin Wang, 
Department of Applied Mathematics, Zhejiang University of Technology, Hangzhou 310023, China}
\email{{{yanlin90118@163.com}};  yanlinwang.math@gmail.com}

\begin{document}
\let\origmaketitle\maketitle
\def\maketitle{
  \begingroup
  \def\uppercasenonmath##1{} 
  \let\MakeUppercase\relax 
  \origmaketitle
  \endgroup
}

\maketitle

\begin{abstract}
In this paper we study the zero-viscosity limit of 2-D Boussinesq equations with vertical viscosity and zero diffusivity,  which is a nonlinear  system with partial dissipation arising in atmospheric sciences and oceanic circulation. The domain is taken as $\mathbb{R}_+^2$ with Navier-type boundary.  We prove  the nonlinear stability of the approximate solution constructed by boundary layer expansion in conormal Sobolev space.  The expansion order and convergence rates for the inviscid limit  are  also identified in this paper.  Our paper extends a partial zero-dissipation limit result of Boussinesq system with full dissipation  by Chae D. [{\it Adv. Math. 203 (2006), no.2, 497–513}]   in the whole space    to  the case with partial dissipation  and  Navier  boundary in the half plane. 



\end{abstract}


\section{Introduction}
We consider the  following 2-D  Boussinesq equations with only  vertical viscosity and zero diffusivity in  $\mathbb{R}_+^2$:
\begin{equation}\label{eq1}
\begin{cases}
& \partial_t \bold{u}^\varepsilon+\bold{u}^\varp\cdot \nabla \bold{u}^\varp+\nabla p^\varp-\varepsilon^2\partial^2_y \bold{u}^\varp=\theta^\varp e_2,\\
& \partial_t\theta^\varp+\bold{u}^\varp\cdot \nabla \theta^\varp=0,\\
&\nabla\cdot \bold{u}^\varp=0,\\
&\bold{u}^\varp|_{t=0}=\bold{u}^\varp_0(x, y), \ \ \theta^\varp|_{t=0}= \theta^\varp_0(x, y),
\end{cases}
\end{equation}
where  $t\geq 0$ and $(x, y)\in \mathbb{R}\times \mathbb{R}_+$ are time and space variables.  Here $\bold{u}^\varepsilon=(u_1^\varepsilon, u_2^\varepsilon)$ and  $p^\varepsilon$ are, respectively, the velocity and pressure of the fluid.  The scalar function $\theta^\varepsilon$  denotes  the  temperature or the  density. We  use  $e_2=(0,1)$ to represent the unit vector in the vertical direction, and $\varepsilon^2$ to denote the kinematic viscosity. 
The Boussinesq 
equations with anisotropic (full or partial) dissipation    play an important role in the study of  atmospheric and oceanographic flows and Raleigh-B\'enard convection, mathematically and physically.  For more backgrounds of the Boussinesq equations, one can refer to \cite{M1,P, CD,Wen}.

In this paper, the Navier-type (slip) boundary condition of the system \eqref{eq1} is given by 
\begin{align}
u_2^\varp=0,\  \  \partial_y u_1^\varp=\alpha u_1^\varp\ \ \  {\rm on} \ \  \{y=0\},\label{navier}
\end{align}
where $\alpha\in \mathbb{R}$ is used to characterize  the tendency of the fluid to slip on the boundary. Here  the initial data given in \eqref{eq1}  and in the sequel  should  satisfy the compatibility conditions on the  boundary and  the divergence free condition. 

Our purpose  is to survey   the zero-viscosity limit behaviour of the Boussinesq system  with partial viscosity and  zero diffusivity  in a half space  satisfying the Navier boundary condition.    It is extremely challenging   to deal with the loss of partial dissipation and the boundary layer effects, especially, for the bounded / unbounded domain with non-slip boundary. In this paper we set our problem in the half plane with the Navier boundary condition \eqref{navier}  for the first step.  To our best knowledge,  this paper is the first one to  consider the  strong zero-viscosity limit of Boussinesq equations with partial viscosity and boundaries,  which is from  physical consideration and involves  layer effects.

Formally, letting $\varepsilon\to 0$, the partially viscous  Boussinesq system \eqref{eq1} is then reduced to the following zero dissipation system in $\mathbb{R}^2_+:$
\begin{equation}\label{var0}
	\begin{cases}
		& \partial_t \bold{u}^0+ \bold{u}^0\cdot \nabla \bold{u}^0+\nabla p^0=\theta^0 e_2,\\
		& \partial_t\theta^0+\bold{u}^0\cdot \nabla \theta^0=0,\\
		&\nabla\cdot \bold{u}^0=0,\\
		&\bold{u}^0|_{t=0}=\bold{u}_0^0(x,y),\ \ \theta^0|_{t=0}=\theta_0^0(x, y),
	\end{cases}
\end{equation}
which is complemented with  the boundary condition: 
\begin{align}
u_2^0|_{y=0}=0.\label{var0b}
\end{align} We  note that the initial data $(\bold{u}^0_0, \theta^0_0)$ satisfy the compatibility condition
$$\nabla\cdot \bu_0^0=0, \ \  \bu_0^0\cdot \bold{n}=0,\ \ {\rm with  } \ \  \bold{n}=(0, -1).$$

 Setting  the zero-dissipation  system \eqref{var0} into the whole space $\mathbb{R}^2$ or a  domain with boundary, the local in time existence theory of \eqref{var0} has been studied in \cite{CN, CKN} or \cite{CH, CI, HH, EJ}, respectively. From \cite{CI, CN},  one can expect the local existence of \eqref{var0} on a time interval $[0, T]$ for initial data satisfying   $\nabla \cdot \bold{u}_0^0=0$  and $(\bold{u}^0_0, \theta_0^0)\in H^{s}(\mathbb{R}^2_+)$  for some $s>2.$  In particular, the local in time existence of \eqref{var0} in $\mathbb{R}_+^2$ holds for 
the  $C^\infty$-smooth initial data \cite{EJ} and for the $C_c^{1,\alpha}(\mathbb{R}^2_+)$ initial data \cite{CH} with $\alpha<\alpha_0$ for some $\alpha_0<1$, both developing singularities in finite time. 
 It is suggested that there are some strong analogues between the 2-D  zero-dissipation   system \eqref{var0} and the 3-D Euler equations (see \cite{CMT, M}).  Therefore, to understand the local existence theory of  the 2-D Boussinesq system with no dissipation \eqref{var0},  one can also refer to  the local existence theory and finite time blow-up criteria  for 3-D incompressible Euler equations \cite{BKM, CC, CH, CHKLSY, Elgindi, ElJ,  Kato, CHH, LuoHou, TaoT}.


The system \eqref{eq1} is a particular  case  of the 2-D Boussinesq system with  anisotropic dissipation, which reads
\begin{equation}\label{eq-aniso}
\begin{cases}
& \partial_t \bold{u}^\varepsilon+\bold{u}^\varp\cdot \nabla \bold{u}^\varp+\nabla p^\varp=\nu_1\partial^2_x \bold{u}^\varp+\nu_2\partial^2_y \bold{u}^\varp+\theta^\varp e_2,\\
& \partial_t\theta^\varp+\bold{u}^\varp\cdot \nabla \theta^\varp=\kappa_1\px^2 \theta^\varp+\kappa_2\py^2 \theta^\varp,\\
&\nabla\cdot \bold{u}^\varp=0,
\end{cases}
\end{equation}
where $\nu_1, \nu_2$ and $\kappa_1, \kappa_2$ are nonnegative constants to characterize viscosity and diffusivity, respectively.  
Recently, extensive progress  has been made on the global  existence theory of the Boussinesq system \eqref{eq-aniso}  with full dissipation or partial  dissipation 
in $\mathbb{R}^2$ (cf. \cite{AH, Aetal,CaoW,CDh, DP, DP1, DWXZ, ES, HKR, HKR2, LWZ, LT, MZ, WX} and the references therein), or in an appropriate  domain with boundary (cf. \cite{DWZZ, HWW, LPZ, LWZ, SZh, ZK} and the references therein).  However, the global existence {theories} for the Boussinesq system with $\nu_2>0$ and $\nu_1=\kappa_1=\kappa_2=0$ (i.e. \eqref{eq1}), or $\nu_1=\nu_2=\kappa_1=\kappa_2=0$ (i.e. \eqref{var0}),  are still  challenging open problems even in the whole space $\mathbb{R}^2$.  For the coupled  advective scalar equation of  $\theta$  in $\eqref{eq-aniso}$ and the related quasigeostrophic equations,   they are also  attractive research topics in PDE theory  \cite{CP, CD,CMT, CC, Ju, RS, W}.

In addition, the zero-dissipation limit behaviour of Boussinesq equations with boundary is a meaningful physical problem, which is associated with the boundary layer theory. Recently,  Jiang et al \cite{JJZ} have studied the zero-diffusivity limit of the 2-D Boussinesq equations  with full dissipation   in the half plane. They have justified the  zero-diffusivity  limit of the velocity $\bu$ and  the  temperature $\theta$, respectively,  in $H^1$ and  $L^2$ norm uniformly on $[0,T]$ with convergence rate.  In \cite{WXie}, Wang and Xie have investigated the $L^2$  zero-dissipation limit of Boussinesq equations in a bounded domain with Navier type boundary for the  velocity and Neumann boundary for the temperature 
using boundary layer expansion. One can  also consult \cite{WW} about the well-posedness of the boundary layer  equation for  a geophysical model or  \cite{FXWW, LXZ}  on the  zero-diffusivity  limit  for the Boussinesq  system  in the weak sense.

The vanishing viscosity problem is an attractive topic in fluid dynamics, due to its physical importance.    It is well known that the famous Prandtl boundary layer equations can be derived from the vanishing viscosity  limit problem of classical incompressible  Navier-Stokes equations with non-slip boundary. The Prandtl equations are  regarded as a leading order of the inner layer approximation.  Similarly, the Prandtl equations can also be derived from the inviscid limit problem of Boussinesq equations with non-slip boundary condition on velocity.   In past decades, much attention have been paid to the well-posedness theory  of Prandtl equations in different function spaces, mainly overcoming the difficulties from the lack of dissipation in horizontal direction.  Based on the well-posedness theory of Prandtl equations, a natural question is that whether the Prandtl layer approximation is stable or not?  In this direction, a lot of  progresses have been made for Navier-Stokes equations ( For example, see \cite{KVW, MY, SC1, SC2, WWZ} and the references therein). It is worthy to investigate the similar questions in ocean dynamic describing by anisotropic Boussinesq equations.  Especially, when the horizontal viscosity is laking, just like the model presented in this paper,  the $L^\infty$ zero-viscosity limit problem is  extremely challenging  in a general domain with non-slip boundary, due to the strong Prandtl layer and the less dissipation.   As a first step, in this paper we set our inviscid limit problem in the half plane with vertical viscosity and  a (slip) Navier boundary.  Meanwhile, our  results can be extend to other anisotropic Boussinesq models  with weak boundary layers in a half space. We hope to  come back to the inviscid limit problem of anisotropic Boussinesq equations in a general domain in the future study.

In this paper we concentrate  on investigating the strong $L^\infty$ zero-viscosity limit of the 
Boussinesq equations \eqref{eq1} with Navier boundary in $\mathbb{R}_+^2.$  To begin with, we construct an approximate solution with outer layer profile away from the boundary  and inner layer profile near the boundary, using boundary layer expansion method.  Due to the structure of the Navier boundary condition, we can derive an approximate solution  with leading profile 
 $$(\bu^0, \theta^0)+\varp(U, \Theta)(t, x, \frac{y}{\varp})+\cdots$$ 
 to characterize the singularity near the boundary as $\varp$ goes to zero. Here we can observe that the $H^2$ norm of $\bold{u}^\varp$ is not uniformly bounded.  The study of zero-viscosity limit is then reduced to the stability analysis of the approximate solution and derive the uniform estimates on $[0, T]$ with $T$ independent of $\varp.$  If the expansion order is larger than or equal to two,  we  can prove the linear stability  in the conormal Sobolev setting by taking advantage of the anisotropic Sobolev embedding inequality, since   the first order derivatives of the inner layer profile are bounded.  In the nonlinear stability analysis, we  introduce an energy function for the perturbation equation around the approximate scheme, involving $L^ {\infty}$-norm of $\theta$ and $\bu$ and their derivatives, which is inspired by the remarkable work \cite{MR}.   This will overcome the difficulties  from  nonlinearity with the help of a precise $L^\infty$ estimates.  More precisely, we deal with  the  ${L^\infty}$ estimates  of  the velocity $\bu$  in spirit of the  $L^\infty$  estimates  for incompressible  Navier-Stokes equations in \cite{MR}, taking advantage of  estimates on Green's function for the operator of  an approximate equation. However,  to close the essential nonlinear estimates, a higher order conormal energy estimate is needed and the order of singularity will increase, compared with the linear argument.  Here we  can establish a uniform convergence estimate through improving  the order of boundary layer expansion. Notice that the  lack of horizontal dissipation can be handled  owing to the divergence free condition and the boundary setting.  Throughout our stability analysis, we can construct  an approximate solution near the boundary with high accuracy.   Different from the inviscid limit problem of Navier-Stokes equations with Navier boundary considered in \cite{MR}, the Boussinesq equations in this paper  are coupling with  an advective scalar equation on $\theta$.   As for the ${L^\infty}$ estimates of  the density (or temperature) $\theta,$  we apply  the maximum principle for a particular  advective scalar equation ( cf. \cite{RS}).  
Here we shall notice that in \cite{MR} the authors have proved the zero-viscosity limit of incompressible Navier-Stokes equations with Navier boundary using compactness argument,  in which the convergence rate  can not be deduced. 
Owing to the boundary layer expansion and stability analysis in this paper,  we can deduce the convergence rate $O(\varp)$ and identify the  expansion order  $K$ for the boundary layer expansion of the partially dissipative Boussinesq system.

Now we state our main result in the following theorem.
\begin{Theorem}\label{maintheorem}
Let $(\bu^0, \theta^0)$ be a solution to \eqref{var0} \eqref{var0b} defined on $[0, T].$
Then 
 there exists
$(\bu^\varp, \theta^\varp)$ a solution to \eqref{eq1} \eqref{navier} defined on $[0, T_1]$ for some $T_1$ independent of $\varp$ and $T_1\leq T$ such that 
\begin{align}
&\sup_{[0, T_1]}(\norm{\bu^\varp-\bu^0}_{L^2(\mathbb{R}_+^2)}+\norm{\theta^\varp-\theta^0}_{L^2(\mathbb{R}_+^2)})\rightarrow 0,\\
&\sup_{[0, T_1]}(\norm{\bu^\varp-\bu^0}_{L^\infty(\mathbb{R}_+^2)}+\norm{\theta^\varp-\theta^0}_{L^\infty(\mathbb{R}_+^2)})\rightarrow 0,
\end{align}
as $\varp$ goes to $0.$
Furthermore, the rate of convergence is $O(\varp).$

\end{Theorem}

\begin{Remark}
The convergence results in Theorem \ref{maintheorem} also hold for Boussinesq system \eqref{eq-aniso} in $\mathbb{R}_+^2$ with either one of the following two conditions:
\begin{enumerate}
\item $\nu_2=\varp^2,\nu_1=0 \ {\rm or}\ \varp^2,\  \kappa_1=0\ {\rm or}\ \varp^2, \kappa_2=0$ and Navier boundary condition \eqref{navier};
\item $\nu_1=\varp^2, \kappa_1=0 \ {\rm or}\ \varp^2,\ \nu_2=\kappa_2=0$ and  boundary condition $\bu \cdot \bold{n}=0$ on $\{y=0\}$.
\end{enumerate}

In addition, this result can be  extended to the inviscid limit problem of the parallel three-dimensional case. 
\end{Remark}

The arrangement of the remaining sections is as following. We devote Section 2 to constructing an approximate solution using boundary layer expansion. We will give the stability analysis of the linearized system around the approximate solution  in conormal Sobolev space in Section 3. In Section 4, we will focus on dealing with the nonlinear terms and deriving  the $L^\infty$ estimates to close the uniform estimates. 

\section{Boundary Layer Expansion}
We construct an approximate solution of system \eqref{eq1} in  the following  form
\begin{equation}\label{app}
	(\bold{u}_a,\theta_a,p_a)=\sum_{i=0}^K\varepsilon^i(\bold{u}^i,\theta^i,p^i)(t,x,y)+\sum_{i=0}^K\varepsilon^i\left(U^i,\Theta^i,P^i\right)\left(t,x,\frac{y}{\varepsilon}\right),
\end{equation}
where $K$ is an arbitrarily large integer, and $U^i:=(U_1^i, U_2^i)$. We use $(\bold{u}^i, \theta^i, p^i)(t, x, y)$ in the first sum to approximate  the outer layer, meanwhile, $\left(U^i,\Theta^i,P^i\right)\left(t,x,\frac{y}{\varepsilon}\right)$ to characterize the inner layer behaviour near the boundary in the second sum. In the following presentation, we denote the fast variable as $z=y/\varp$ for simplicity.

We shall impose the fast decay condition as following          
\begin{align}
\left(U^i,\Theta^i,P^i\right)(t, x, z)\rightarrow 0, \ \ {\rm as}\ \ z\rightarrow +\infty.\label{fastdec}
\end{align}
In addition, the matched boundary conditions for Navier boundary \eqref{navier} read
\begin{equation}\label{matching}
\begin{split}
&u_2^i(t, x, 0)+U^i_2(t, x, 0)=0,\ {\rm for\   all } \ i\geq 0,\\
&\partial_z U_1^0(t,x,0)=0,\\
&\partial_yu_1^i(t, x, 0)+\partial_z U_1^{i+1}(t, x, 0)=\alpha(u_1^i+U^i_1)(t, x, 0), \ i\geq 0.
\end{split}
\end{equation}
In the sequel, we denote $\Gamma f$ as 
$$\Gamma f=f(t, x, y)|_{y=0}.$$

We expect the leading order of the outer layer in the approximation \eqref{app} to be the solution of the zero viscosity Boussinesq system \eqref{var0} \eqref{var0b}, whose local well-posedness theory can be deduced from \cite{CN, CI} for initial data satisfying   $\nabla \cdot \bold{u}_0^0=0$  and $(\bold{u}^0_0, \theta_0^0)\in H^{s}(\mathbb{R}^2_+)$  for some $s>2.$ In fact, substituting the approximate scheme \eqref{app} into the viscous Boussinesq system
\eqref{eq1} and collecting the $O(1)$ order terms, consequently we can get the leading order terms for the outer layer satisfying \eqref{var0} by taking $z\rightarrow +\infty.$

In the inner zone, by collecting the $O(\varp^{-1})$ terms, one has
\begin{align}
&(\Gamma u_2^0+U_2^0)\partial_z U_1^0=0,\notag\\
&(\Gamma u_2^0+U_2^0)\partial_z U_2^0+\partial_z P^0=0,\notag\\
&\partial_z U_2^0=0,\notag\\
&(\Gamma u_2^0+U_2^0)\partial_z \Theta^0=0.\notag
\end{align}
Then we can deduce from \eqref{fastdec} and $\eqref{matching}_1$ that
\begin{align}
U_2^0(t, x, z)=0,\ \ P^0(t, x, z)=0.
\end{align}
In turn, $\Gamma u_2^0=0$ holds.

Now, we collect the $O(1)$ terms for the first equation of velocity in the inner zone to obtain
\begin{align}
&\partial_t(\Gamma u_1^0+U^0_1)+(\Gamma u_1^0+U_1^0)\partial_x(\Ga u_1^0+U_1^0 )\notag\\
&+(\Gamma u_2^1+U_2^1+z\Ga \py u_2^0)\pz (\Ga u_1^0+U_1^0 )
+\px (\Ga p^0)-\partial_{zz}(\Ga u_1^0+U_1^0)=0,\label{nov-1}
\end{align} 
Combining with the $O(1)$ terms from the divergence free condition in the inner layer,
\begin{align}
&\px (\Ga u_1^0+U_1^0)+\pz(\Gamma u_2^1+U_2^1+z\Ga \py u_2^0) =0,\label{nov-2}
\end{align}
one has a closed system for $\Gamma u_1^0+U^0_1$ and 
$\Gamma u_2^1+U_2^1+z\Ga \py u_2^0$, together with the boundary condition
\begin{align}
\pz U_1^0(t,x,0)=0.\label{U1b}
\end{align}
We notice that the equations \eqref{nov-1} \eqref{nov-2} and \eqref{U1b} have a trivial solution
$$U_1^0=0.$$
 In the following analysis, we shall assume that the approximate solutions \eqref{app}  satisfy $U_1^0=0.$
 
  Clearly, from the $O(1)$ order terms of the  second velocity equation and the transport equation, respectively,
\begin{align}
&\partial_t(\Gamma u_2^0+U_2^0)+(\Gamma u_1^0+U_1^0)\px (\Gamma u_2^0+U_2^0)\notag\\
&+(\Gamma u_2^0+U_2^0)\partial_z(\Gamma u_2^1+U_2^1+z\Gamma \py u_2^0)+\Gamma \py p^0+\pz P^1=\Gamma \theta^0+\Theta^0,\notag\\
&\partial_t(\Gamma \theta^0+\Theta^0)+(\Gamma u_1^0+U_1^0)\px (\Gamma \theta^0+\Theta^0)+(\Gamma u_2^1+U_2^1+z\Gamma \py u_2^0)\pz \Theta^0=0,\notag
\end{align}
 the leading order for the temperature  (or the density) in the inner zone  has a trivial solution 
$$\Theta^0=0.$$
Consequently,  we also have $P^1=0.$

\begin{Remark}
In the study of zero-viscosity limit for incompressible Navier-Stokes equations with non-slip boundary, one derive the system \eqref{nov-1} \eqref{nov-2} with boundary condition
\begin{align}\Ga u_1^0+U_1^0(t,x, 0)=0,\label{nonslip}\end{align}
which just forms the famous Prandtl system. For the local well-posedness theory  of Prandtl system\eqref{nov-1} \eqref{nov-2} \eqref{nonslip} or the  inviscid limit of Navier-Stokes equations with non-slip boundary, one may refer to \cite{DG, KVW, LMY, MY, SC1, SC2, WWZ} and the references therein.
\end{Remark}

Till now, we get the leading order for both inner zone and outer zone for the approximate solution. i.e. $(\bold{u}^0, p^0, \theta^0) \ {\rm and} \ (U^0, P^0, \Theta^0)$ are solved. {{Next, following the standard 
method on constructing approximate solution by boundary layer expansion (see \cite{GHR, WXie, Wang} and the references therein for example), one can obtain the terms $(\bold{u}^i, p^i, \theta^i)$ and $(U^i, P^i, \Theta^i)$ order by order. Indeed,  we collect the $O(\varp^i), i\geq1$ terms in the outer layer expansion  to obtain
\begin{align*}
&\pt \bu^i+\bu^i\cdot\nabla \bu^0+\bu^0\cdot \nabla \bu^i+\nabla p^i=\theta^i e_2+\py^2 \bu^{i-2}+f_u^i,\\
&\pt \theta^i+\bu^i\cdot\nabla\theta^0+\bu^0\cdot\nabla\theta^i=f_\theta^i,\\
&\nabla\cdot\bu^i=0,
\end{align*}
where $\bu^{-1}=0, $ and $f_u^i, f_\theta^i$ are depending on $\bu^j, \theta^j$ with $j\leq i-1.$ This, together with the boundary condition $u_2^i|_{y=0}=0,$  forms a linear closed system of $(\bu^i, \theta^i, p^i).$  Hence  the outer layer terms have been determined.   For the 
inner layer expansion in the order $O(\varp^i), i\geq1,$  we have that
\begin{align*}
&\pt (\Gamma u_1^i+U_1^i)+(\Gamma u_1^i+U_1^i)\px (\Gamma u^0_1+U_1^0)+ (\Gamma u^0_1+U_1^0)\px(\Gamma u_1^i+U_1^i)\\
& \ \ \ +(\Gamma u_2^{i+1}+U_2^{i+1}+z\Gamma \py u_2^i)\pz(\Gamma u_1^0+U_1^0)+\Gamma \px p^i+\px P^i=\partial_{zz}(\Gamma u_1^i+U_1^i)+F_{U_1}^i,\\
&\pt(\Gamma u_2^i+ U_2^i)
+ (\Gamma u^0_1+U_1^0)\px(\Gamma u_2^i+U_2^i)+\Gamma \py p^i+\pz P^{i+1}\\
&\ \ \ +(\Gamma u_2^1+z\Gamma\py u_2^1+U_2^1)\pz (\Gamma u_2^i+U_2^i+z\Gamma \py u_2^i)\\
&\ \ \ +(\Gamma u_2^i+U_2^i+z\Gamma \py u_2^{i-1})\pz(\Gamma u_2^1+z\Gamma\py u_2^1+U_2^1) \\
&=\partial_{zz}(\Gamma u^i_2+U_2^i)+\Theta^i +\Gamma \theta^i+F_{U_2}^i,\\
&\pt(\Gamma \theta^i+\Theta^i)+(\Gamma u_1^i+U_1^i)\px(\Gamma \theta^0+\Theta^0)+(\Gamma u_1^0+U^0)\px(\Gamma\theta^i+\Theta^i)\\
&\ \ \ +(\Gamma u_2^1+U_2^1+z\Gamma \py u_2^0)\pz(\Gamma \theta^i+\Theta^i+z\Gamma\py \theta^{i-1})\\
&\ \ \ +(\Gamma u_2^i+U_2^i+z\Gamma \py u_2^{i-1})\pz(\Gamma \theta^1+\Theta^1+z\Gamma\py\theta)=F_{\Theta}^i,\\
&(\Gamma u_2^{i+1}+U_2^{i+1}+z\Gamma \py u_2^i)=-\int_0^z \px(\Gamma u_1^i+U_1^i)(t, x, \tau)d\tau ,
\end{align*}
where $F_{U_1}^i, F_{U_2}^i,F_{\Theta}^i,$ are depending on $U^j, \Theta^j, \bu^j, \theta^j$ with $j\leq i-1.$ Here the last equality is from the divergence free condition. Recall that $P^1=0$ and the boundary conditions for the inner layer satisfy \eqref{matching}. Then, order by order,  we obtain a linear  closed differential and  integral   system 
 for $(U^i, \Theta^i, P^{i+1}), i\geq 1.$ }}
Now  we  conclude that the approximate solutions $(\bold{u}_a, p_a, \theta_a)$ for \eqref{eq1} can be expressed  in the form
\begin{align}
(\bold{u}^0, p^0, \theta^0)(t,x,y)+\varp\left((\bold{u}^1, p^1, \theta^1)(t,x,y)+(U^1, P^1, \Theta^1)(t, x, \frac{y}{\varp})\right)+\cdots.\notag
\end{align}

Let us denote 
$$\omega^0=\py u_1^0-\px u_2^0, \ \ \ {\rm and}\ \ \ \omega_0^0=\omega^0(t=0, x, y).$$
Then we can gather our results on boundary layer approximations of the solutions for the system 
\eqref{eq1} and \eqref{navier} in the following theorem.
\begin{Theorem}\label{appsolu}
Let $K\in \mathbb{N}_+.$ For any initial data $(\theta_0^0, \bold{u}_0^0)$ satisfying $\nabla\cdot \bold{u}_0^0=0,$ $(\bu^0_0,  \theta_0^0)\in H^{s}(\mathbb{R}^2_{+})$  for some $s>2$ and  some compatibility conditions on $\{y=0\}$,  
then there exists $T>0$
and a smooth approximate solution $(\bold{u}_a, \theta_a)$ of \eqref{eq1} under the form
\eqref{app} with order $K$ such that
\begin{itemize}
\item [i)]
we have $(\bold{u}^0, \theta^0)\in C^0([0,T], H^{s}(\mathbb{R}^2_+))$
 as a solution of the inviscid  system \eqref{var0} with initial data $(\bold{u}_0^0, \theta_0^0);$
\item [ii)]
for all $1\leq i \leq K,$ $(\bold{u}^i, \theta^i)\in C^0([0, T], H^{s}(\mathbb{R}^2_+));$ 
\item [iii)]
for all $0\leq i\leq K,$ $(U^i, P^i, \Theta^i)$ are solved at least locally and satisfy the fast decay property \eqref{fastdec}  with respect to the last variable.
\item [iv)] we consider $(\bu^\varp, \theta^\varp)$ a solution to \eqref{eq1}, and 
denote the error terms $(\bu, p, \theta)$ as following
\begin{align}
\bu=(u_1, u_2)=\bu^\varp-\bu_a,\ \ \ p=p^\varp -p_a, \ \ \ \theta=\theta^\varp-\theta_a,\notag
\end{align}
where $\bu_{a}$ naturally  satisfies $$\nabla \cdot \bu_a=0,\ \ \  \bu_{a}\cdot \bold{n}|_{y=0}=0,$$ with $\bold{n}=(0, -1).$

Then $(\bu, \theta)$ satisfies the system of equations
\begin{subequations}\label{erroreq}
\begin{align}
&\pt \bu+\bu\cdot\nabla(\bu+\bu_a)+ \bu_a\cdot\nabla \bu +\nabla p- \varp^2\py^2\bu
=\theta e_2+\varp^K R_{\bu},\label{err-a}\\
&\pt \theta+\bu \cdot \nabla (\theta+\theta_a)+\bu_a\cdot\nabla \theta=\varp^K R_\theta, \label{err-b}\\
&\nabla\cdot \bu =0,\label{err-c}\\
&\py u_1=\alpha u_1, \ \ \ u_2=0, \ \ \ {\rm on} \ \ \  \{y=0\},
\end{align}
\end{subequations}
where $R_{\bu}:=(R_{u_1}, R_{u_2})$ and  $R_\theta$ are remainders satisfying
\begin{align}
\sup_{[0, T]}\norm{\nabla ^\beta R_{\bu, \theta}}\leq C_a \varp^{-\beta_2}, \ \ \ \forall \beta=(\beta_1, \beta_2)\in \mathbb{N}^2,
\end{align}
with $C_a>0$ independent of $\varp.$ {{ Here and thereafter we use the symbol   $\|\cdot\|$ to represent  the standard $L^2$ norm $\|\cdot\|_{L^2(\mathbb{R}^2_+)}$. }}

\end{itemize}

\end{Theorem}

\section{Linear Stability Estimates}
This section is devoted to the stability analysis of the approximate solution constructed by boundary layer expansion in Theorem \ref{appsolu}. Due to the essential challenge in nonlinear energy estimates  caused by the boundary layer effects  and the loss of dissipation in $x$-direction, we first consider the linear stability of the boundary layer approximation. 
In the following writing, we use the notation 
$$\Omega:=\mathbb{R}_+^2.$$

To begin with, we  linearize  the system \eqref{erroreq}  around the approximate solution to obtain that,  in $\mathbb{R}_+^2,$
\begin{subequations}\label{erroreq-linear}
\begin{align}
&\pt \bu+\bu\cdot\nabla \bu_a + \bu_a\cdot\nabla \bu +\nabla p- \varp^2\py^2\bu
=\theta e_2+\varp^K R_{\bu},\label{errl-a}\\
&\pt \theta+\bu \cdot \nabla \theta_a +\bu_a\cdot\nabla \theta=\varp^K R_\theta, \label{errl-b}\\
&\nabla\cdot \bu =0,\label{errl-c}\\
&\py u_1=\alpha u_1, \ \ \ u_2=0,  \ \ \ {\rm on} \ \ \  \{y=0\},\label{errl-d}
\end{align}
\end{subequations}
where $R_{\bu}:=(R_{u_1}, R_{u_2})$ and $R_\theta$ are remainders satisfying
\begin{align}
\sup_{[0, T]}\norm{\nabla ^\beta R_{\bu, \theta}}\leq C_a \varp^{-\beta_2}, \ \ \ \forall \beta=(\beta_1, \beta_2)\in \mathbb{N}^2,
\end{align}
with $C_a>0$ independent of $\varp.$  The initial data of the system \eqref{erroreq}  and \eqref{erroreq-linear} are given by 
\begin{align}
\bu|_{t=0}=\varp^{K+1} \bu_0, \ \ \ \theta|_{t=0}=\varp^{K+1} \theta_0.\label{ini-linear}
\end{align}
{{Here we can always write  the initial data $(\bold{u}_0^\varepsilon, \theta_0^\varepsilon)$ of the system \eqref{eq1} as 
$$(\bold{u}_0^\varepsilon, \theta_0^\varepsilon)=\left(\sum_{i=0}^K \varepsilon^i(\bold{u}^i_0, \theta^i_0)\right)+\varepsilon^{K+1}(\bold{u}_0, \theta_0),$$
where $(\bold{u}^i_0, \theta^i_0)$ is taken  as  the initial data for each approximate profile  in \eqref{app}. Hence $\varepsilon^{K+1}(\bold{u}_0, \theta_0)$ is regarded as 
the initial data of the perturbation equations \eqref{erroreq} or \eqref{erroreq-linear}.}}

In the remaining parts of this section, we will give the stability analysis of the linear system \eqref{erroreq-linear} for  $\alpha \in \mathbb{R}.$

\subsection{Preliminaries}
Before we process the uniform estimates, some basic notations and useful inequalities will be given in this subsection. We will use $\lesssim$ to denote 
$\leq C(\cdot)$ or $\leq C_a(\cdot)$ for a generic constant $C$ or a constant $C_a$ depending on the approximate solution, but both independent of $\varp$.
The standard Sobolev norm is denoted by $\norm{\cdot}_{s}$ for  $\norm{\cdot}_{H^s}$  with $s\geq 0.$ In particular, $\norm{\cdot}$ is for $\norm{\cdot}_{L^2}$.

We introduce conormal operator $Z^{\beta}=Z_1^{\beta_1}Z_2^{\beta_2}$ with $\beta=( \beta_1, \beta_2)\in \mathbb{N}^2,$
where 
\begin{align}
 Z_1=\px, \ \ \ \ Z_2=\varphi(y)\py,
\end{align}
and $\varphi(y)$ is a smooth function satisfying $\varphi(0)=0$ and $\varphi'(0)>0$, such as
\begin{eqnarray}
\varphi(y)=\frac{y}{y+1}.
\end{eqnarray}
In the sequel, we denote 
$$Z^k=Z^{\beta},\ \ \ {\rm with}\ \ |\beta|=k, \ k\in \mathbb{N}.$$
The conormal Sobolev space $H_{co}^s$ for $s\in \mathbb{N}$ is defined by
\begin{align}
{{
H_{co}^s:=\left\{u\in L^2_{x,y}(\mathbb{R}^2_+):\norm{u}_{H_{co}^s}^2=\sum_{|\beta|\leq s}\norm{Z_1^{\beta_1}Z_2^{\beta_2} u}^2_{L^2_{x,y}(\mathbb{R}_{+}^2)}<\infty \right\}.  }} \label{conormal}
\end{align}
In spirit of  this setting, we denote the $H_{co}^s$-inner product by $\langle \cdot , \cdot\rangle_{H_{co}^s}.$
Set $$\norm{u}_{k,\infty}=\sum_{|\beta|\leq k}\norm{Z^{\beta} u}_{L^\infty},$$
and we say that $u\in W_{co}^{k,\infty}$ if $\norm{u}_{k, \infty}$ is finite.

Now we state the anisotropic Sobolev embedding inequality and some useful estimates  (cf. \cite{MR}) in the following lemmas.
\begin{Lemma}\label{le1}
For $m_0\geq 1, m_0\in \mathbb{N},$ and conormal Sobolev norm defined in \eqref{conormal},  we have 
\begin{align}
\norm{u}^2_{L^{\infty}}\les  \norm{\py u}_{H_{co}^{m_0}} \norm{ u}_{H_{co}^{m_0}}+ \norm{ u}^2_{H_{co}^{m_0}}. \label{coinfty}
\end{align}

\end{Lemma}

Let us denote the vorticity $\omega$ as 
$$\omega =\textup{curl}\  {\bu}=\py u_1-\px u_2,$$ 
and  set 
\begin{align}\eta=\omega-\alpha u_1.\label{etadef}
\end{align}

\begin{Lemma}[cf.  \cite{MR}, Proposition 12]\label{le2}
In $\mathbb{R}^2_+$,  for $m_0\geq 1,$ we have
\begin{align}
&\norm{\bu}_{W^{1,\infty}}\les \norm{\bu}_{H^{m_0+2}_{co}}+\norm{\eta}_{H^{m_0+1}_{co}}+\norm{\eta}_{L^\infty}, \label{etainfty-1}\\
&\norm{\bu}_{2,\infty}\les \norm{\bu}_{H_{co}^{m_0+3}}+\norm{\eta}_{H^{m_0+2}_{co}},\label{etainfty-2}\\
&\norm{\nabla \bu}_{1,\infty}\les \norm{\bu}_{H_{co}^{m_0+3}}+\norm{\eta}_{H^{m_0+3}_{co}}+\norm{\eta}_{1,\infty}.\label{etainfty-3}
\end{align}
\end{Lemma}

\subsection{Uniform Estimates for $\alpha\in \mathbb{R}$.}  In this subsection, we will derive  the uniform estimates in conormal Sobolev  space.  To begin with, we state our main results  for  the  linear system \eqref{erroreq-linear} with partial viscosity.
\begin{Theorem}\label{thmlinearco}
Let $m\geq 1,$ and $(\theta^0_0, \bu^0_0)\in C^\infty(\mathbb{R}^2_+)$ be  initial data for \eqref{var0} satisfying $\nabla\cdot \bu_0^0$ and some compatibility conditions on $\{y=0\}$. Let $K\in \mathbb{N}_+, K> m$ and $(\bu_a, p_a, \theta_a)$ an approximate solution at order $K$ given by Theorem \ref{appsolu}. Then,  for every $\varp \in (0, 1),$  there exists a $T_0>0$ such that  the solution of \eqref{erroreq-linear}--\eqref{ini-linear} defined on $[0, T_0]$   satisfies the estimate
\begin{align}
\norm{(\bu, \theta)}_{\hcom}+\norm{\py(\bu, \theta)}_{\hcoml}\les \varp^{K-m}.\label{thm-co}
\end{align}
Therefore, it holds that 
\begin{align}
\norm{(\bu, \theta)}_{L^\infty(\mathbb{R}^2_+)}\les \varp^{K-1}.\label{thm-co1}
\end{align}

\end{Theorem}  

The proof of the Theorem \ref{thmlinearco} can be deduced from the following $L^2$-estimates, conormal energy estimates, normal derivatives estimates and pressure estimates.

{\it \textbf{$L^2$-estimates.}}
\begin{Lemma}\label{le-l2}
Let $(\bu_a, \theta_a)$ be the approximate solution in  Theorem \ref{appsolu}. Given $\alpha\in \mathbb{R},$    if $(\bu, \theta)$ is the solution of  the linear system \eqref{erroreq-linear} defined on $[0, T],$   then it holds that
\begin{align}\label{l2est}
\norm{(\bu, \theta)(t, x, y)}^2+c_0\varp^2  \int_0^T \norm{\py u(\tau, x, y)}^2 d\tau \les  \varp^{2K}
\end{align}
for some positive constant $c_0.$
\end{Lemma}

\begin{Remark}
The $L^2$ estimates \eqref{l2est} also holds for the nonlinear system \eqref{erroreq},
due to $\nabla \cdot \bu=0$,  and $\bu\cdot \bold{n}|_{y=0}=0.$ 
\end{Remark}

\begin{proof}
We multiply the equation \eqref{errl-a} and \eqref{errl-b}, respectively,  by $\bu$
and $\theta.$ Then,  adding the resulting equations together and integrating over $\Omega,$  together with integration by parts, we can obtain that 
\begin{align}
&\frac{d}{dt}\int_{\Omega} (\bu^2+\theta^2)+\varp^2 \int_{\Omega}|\py \bu|^2
+\alpha \varp^2 \int_{\partial \Omega} u_1^2\notag\\
&=-\int_{\Omega}( \bu\cdot\nabla \bu_a \cdot \bu+\bu \cdot \nabla \theta_a \theta)+{{\int_{\Omega} \theta e_2\cdot\bold{u}}}+\varp^K \int_{\Omega} (R_{\bu}\bu+R_\theta \theta),
\end{align}
since 
\begin{align*}
& \int_{\Omega}\nabla p \bu=-\int_{\Omega} \nabla\cdot \bu p+\int_{\partial\Omega} p \bu\cdot \bold{n}=0,\\
&\int_{\Omega} \bu_{a}\cdot \nabla \bu \bu=\int_{\Omega} \bu_{a} \cdot \nabla |\bu|^2 =-\int_{\Omega} \nabla \cdot \bu_{a} |\bu|^2+\int_{\partial \Omega} \bu_{a}\cdot \bold{n} |\bu|^2=0,\\
&\int_{\Omega}\py^2\bu \bu=-\int_{\Omega}|\py \bu|^2 -\int_{\partial \Omega} \py u_1 u_1=-\int_{\Omega} |\py \bu|^2-\alpha\int_{\partial \Omega}u_1^2,\\
&\int_{\Omega}\bu_{a}\cdot \nabla \theta \theta=-\int_{\Omega}\nabla\cdot \bu_a \theta^2+\int_{\partial \Omega} \bu_a \cdot \bold{n} \theta^2=0.
\end{align*}

Due to 
$$\int_{\Omega}( \bu\cdot\nabla \bu_a \cdot \bu+\bu \cdot \nabla \theta_a \theta) \leq \norm{\nabla \bu_a}_{L^\infty}\int_{\Omega}|\bu|^2+\norm{\nabla \theta_a}_{L^\infty}\int_{\Omega}|\bu||\theta|\leq C_a \int_{\Omega}(|\bu|^2+|\theta|^2),$$
and $${\int_{\Omega}\theta e_2\cdot \bold{u}\leq\int_{\Omega}(|\bu|^2+|\theta|^2), }$$
then we have, for $\alpha\geq 0,$
\begin{align}
&\frac{d}{dt}\int_{\Omega} (\bu^2+\theta^2)+\varp^2 \int_{\Omega}|\py \bu|^2
+\alpha \varp^2 \int_{\partial \Omega} u_1^2\notag\\
&\les \int_{\Omega} (\bu^2+\theta^2)+ \varp^{2K},\label{l2est-1}
\end{align}
where the Cauchy inequality has been used.   However, for $\alpha<0,$ the boundary term is not good in the energy estimates. 
Using trace theorem in $\mathbb{R}_+^2,$ we have 
\begin{align}
|\alpha| \varp^2 \int_{\partial \Omega} u_1^2\leq |\alpha| \varp^2 \norm{\py u_1}\norm{u_1}.\notag
\end{align}
Then, together with the Young inequality,  one has
\begin{align}
&\frac{d}{dt}\int_{\Omega} (\bu^2+\theta^2)+c_0 \varp^2 \int_{\Omega}|\py \bu|^2\les \int_{\Omega} (\bu^2+\theta^2)+ \varp^{2K}.\label{l2est-2}
\end{align}
Hence \eqref{l2est} follows from \eqref{l2est-2} and the  Gronwall's inequality.
\end{proof}

{\it \textbf{Conormal energy estimates.} } 
\begin{Lemma}\label{le-co}
Let $(\bu_a, \theta_a)$ be the approximate solution in Theorem \ref{appsolu}. Given $\alpha\in \mathbb{R},$   assume that $(\bu, \theta)$ is a solution of  the linear system \eqref{erroreq-linear} defined on $[0, T],$   then, for $m\geq 1,$  we have
\begin{align}
&\frac{d}{dt}(\norm{(\bu, \theta)}^2_{H_{co}^m})+c_0\varp^2\norm{\py\bu}^2_{\hcom}\notag\\
&\les \norm{\py \bu}^2_{H_{co}^{m-1}}+\norm{\theta}^2_{\hcom}+\norm{\bu}_{\hcom}^2+\norm{\nabla p}_{H_{co}^{m-1}}\norm{\bu}_{\hcom}+{\color{teal}{\varp^{2K-2m}}},\label{coest}
\end{align}
where $c_0$ is some positive constant independent of $\varp.$
\end{Lemma}

\begin{proof}
Acting $Z^\beta$ with $|\beta|\leq m$ on the equations \eqref{errl-a} and \eqref{errl-b}, we obtain that
\begin{align}
&\pt \zb\bu+\bu\cdot\nabla\zb \bu_a+ \bu_a\cdot\nabla\zb \bu +\nabla\zb p- \varp^2\py^2Z^\beta\bu
=\zb\theta e_2+\mathcal{C_\bu}+\varp^K \zb R_{\bu},\label{errc-a}\\
&\pt \zb \theta+ \bu \cdot \nabla \zb \theta_a+\bu_a\cdot\nabla \zb \theta=\mathcal{C}_\theta+\varp^K \zb R_\theta, \label{errc-b}
\end{align}
where the commutators  $\mathcal{C}_\bu, \mathcal{C}_\theta$ are defined as following 
\begin{align}
&\mathcal{C}_\bu=-[Z^\beta, \bu\cdot \nabla]\bu_{a}-[Z^\beta, \bu_a\cdot \nabla] \bu-[Z^\beta, \nabla]p+{\varp^2[Z^\beta, \py^2]\bu}=:\mathcal{C}_1+\mathcal{C}_2+\mathcal{C}_3+\mathcal{C}_4,\notag\\
&\mathcal{C}_\theta=-[Z^\beta, \bu\cdot\nabla]\theta_a-[\zb, \bu_a\cdot \nabla]\theta:=\mathcal{C}_5+\mathcal{C}_6.\notag
\end{align}
Now we calculate $\eqref{errc-a}\times\zb \bu+\eqref{errc-b}\times \zb \theta$ and then  integrate the resulting equality  over $\Omega$ to get 
\begin{align}
&\frac{1}{2}\frac{d}{dt}\int_{\Omega}(|\zb \bu|^2+|\zb \theta|^2)+\varp^2\int_{\Omega} |\py\zb \bu|^2+\varp^2\int_{\partial\Omega}
\py \zb u_1 \zb u_1\notag\\
&=-\int_{\Omega}(\bu\cdot\nabla\zb\bu_a \zb\bu+ \bu\cdot \nabla Z^\beta \theta_a \zb \theta)-\int_{\Omega}\nabla\zb p \zb\bu+\int_{\Omega} \zb \theta e_2\zb \bu\notag\\
&\ \ \ +\int_{\Omega} \mathcal{C_\bu}\zb\bu+
\int_{\Omega} \mathcal{C}_\theta \zb \theta+\varp^K \int_{\Omega}( \zb R_{\bu}\zb\bu+\zb R_\theta \zb \theta):=\sum_{j=1}^6 I_i.\notag
\end{align}
We deduce from Navier boundary condition \eqref{navier} that
\begin{align}
&\int_{\partial \Omega} \py \zb u_1 \zb u_1=\int_{\partial \Omega} \zb \py u_1 \zb u_1+\int_{\partial \Omega} [\py, \zb] u_1 \zb u_1\notag\\
&= \alpha\int_{\partial \Omega} |\zb u_1 |^2=\alpha\int_{\partial \Omega} |\px^{|\beta|} u_1 |^2=\alpha \norm{u_1}^2_{\hcom(\partial \Omega)}, \label{bterm}
\end{align}
due to 
\begin{align}
&[\py, Z_1 ]=0,\notag\\
&\zb u_1|_{\partial \Omega}=0, \ {\rm if} \ \beta_2>0.\notag
\end{align}
Here we note that if  the boundary $\partial \Omega$ is not flat, $Z_2 \bu|_{\partial \Omega}=0$ does not hold. 
Hence, by using trace theorem,  we can derive that 
\begin{align}
|\alpha|\int_{\partial \Omega} |Z^\beta u_1|^2\leq C |\alpha| \norm{\py u_1}_{\hcom}\norm{u_1}_{\hcom}.
\end{align}
Then we estimate $I_i$ $(1\leq i\leq 6)$ one by one. For $I_1$, it holds that 
{{\begin{align}
|I_1|&\les |\nabla \zb \bu_a|_{L^\infty}\int_{\Omega} \bu\zb\bu+|(\zb\px\theta_a, \py\zb\theta_a)|_{L^\infty} \into( u_1\zb\theta+u_2\zb\theta)\notag\\
&\les (|Z^{\beta-1} \py\bu_a|_{L^\infty}+|\zb\py \bu_a |_{L^\infty}+|(\zb\px\theta_a, \py\zb\theta_a)|_{L^\infty} )\int_{\Omega} (|\bu|^2+|\zb\bu|^2 +|\zb\theta|^2)\notag\\
&\les\int_{\Omega} (|\bu|^2+|\zb\bu|^2 +|\zb\theta|^2)\les \norm{(\bu,\theta)}_{\hcom}^2,\label{i1}
\end{align}}}
due to the leading order term  $(U^0_1, U_2^0)$ and $\Theta^0$ vanishing in the inner layer expansion.
Note that
\begin{align}
I_2&=-\int_{\Omega}\nabla\zb p \zb\bu=\int_{\Omega}\zb p \nabla\cdot \zb \bu
-\int_{\partial \Omega}\zb p\zb\bu\cdot \bold{n}=-\int_{\Omega}\zb p ([\zb, \nabla\cdot]\bu),\notag
\end{align}
and 
\begin{align}
&[Z_2, \nabla\cdot] \bu=-\varphi'\py u_2=\varphi'\px u_1,\notag\\
&[Z_1, \nabla\cdot]\bu=0.\notag
\end{align}
Hence we obtain 
\begin{align}
|I_2|\les \int_{\Omega}( |\zb \bu||\zb p|)\les \norm{\bu}_{\hcom}\norm{\nabla p}_{\hcom} \label{i2}.
\end{align}
By virtue of the Cauchy inequality, we have
\begin{align}
&|I_3|\les \into (|\zb \theta|^2+|\zb \bu|^2)\les \norm{\bu}^2_{\hcom}+\norm{\theta}_{\hcom}^2,\label{i3}\\
&|I_6|\les \into(|\zb \theta|^2+|\zb \bu|^2) +\varp^{2K}\les  \norm{\bu}^2_{\hcom}+\norm{\theta}_{\hcom}^2+{{\varp^{2K-2m}.}}\label{i6}
\end{align}
As for the term $I_4$ containing commutators $\mathcal{C}_i$ $(1\leq i\leq 4)$, 
we have 
\begin{align}
\mathcal{C}_1&=-[\zb, \bu\cdot \nabla]\bu_a=-\left(\zb(\bu\cdot\nabla\bu_a)-\bu\cdot\nabla\zb\bu_a\right)\notag\\
&=-\left(\sum_{\gamma+\zeta=\beta, \gamma\neq 0} c^1_{\gamma,\zeta}Z^{\gamma} \bu\cdot Z^{\zeta}\nabla \bu_a+\bu\cdot [\zb, \nabla] \bu_a\right),\notag\\
\mathcal{C}_2&=-[\zb, \bu_a\cdot\nabla]\bu=-\left(\zb(\bu_a\cdot\nabla\bu)-\bu_a\cdot\nabla\zb\bu\right)\notag\\
&=-\left(\sum_{\gamma+\zeta=\beta, \gamma\neq 0} c^2_{\gamma,\zeta}Z^{\gamma} \bu_a\cdot Z^{\zeta}\nabla \bu+\bu_a\cdot [\zb, \nabla] \bu\right),\notag
\end{align}
{{where  $c^1_{\gamma,\zeta}, c^2_{\gamma,\zeta}$ is, respectively,  from the expansion of $\zb(\bu\cdot\nabla\bu_a)$ and $\zb(\bu_a\cdot\nabla\bu)$ by Leibniz formula}}. 
We  obtain
\begin{align*}
&\into \mathcal{C}_1\zb \bu\les \|\bu\|^2_{H^m_{co}},\\ 
&\into \mathcal{C}_2 \zb \bu\les\|\bu\|^2_{H^m_{co}}+\|\py\bu\|^2_{H_{co}^{m-1}},\notag
\end{align*}
due to 
\begin{align}
&|\nabla \bu_a|_{L^\infty}\leq C_a, \ \ \ 
{u_{a2}|_{\partial \Omega}=0} \ \  ({\rm  i.e.} \ \  \Gamma u_{a2}=0), \notag\end{align}
and
\begin{align} \left \Vert u_{a2}\py Z^{|\beta|-1} \bu \right \Vert=\left \Vert \frac{u_{a2}}{\varphi(y)}\varphi(y)\py Z^{|\beta|-1}\bu \right \Vert\leq |u_{a2}|_{W^{1,\infty}}\norm{Z^{|\beta|} \bu}.\notag
\end{align}
One can easily get 
\begin{align}
\norm{\mathcal{C}_3}\les\norm{ Z^{|\beta|-1}\nabla p} \les \norm{\nabla p}_{H^{m-1}_{co}}.
\end{align}
Together with 
\begin{align}
[Z_1, \py^2]\bu=0, \ \ \ {[Z_2, \py^2]\bu =-\varphi''\py \bu-2\varphi'\py^2 \bu},\ \ \ Z_2 \bu |_{\partial \Omega}=0,\notag
\end{align}
{{we can derive the estimates for the most difficult case  when  $$\mathcal{C}_4=[Z_2^\beta, \py^2]\bu=-\varepsilon^2\left(\py(\beta\varphi'\varphi^{\beta-1}\py^\beta  \bu)+\beta\varphi'\varphi^{\beta-1}\py^\beta \py \bu \right),$$ using  integration by parts:
\begin{align}
&\into \mathcal{C}_4 \zb \bu -\varepsilon^2\int_{\Omega}(\beta\varphi'Z^{\beta-1}\py \bu)^2\notag\\
&=\varepsilon^2 \int (\beta \varphi'Z^{\beta-1}\py \bu)\cdot (\zb\py\bu)-\varepsilon^2 \int_{\Omega}(\beta \varphi'\varphi^{\beta-1}\py^\beta \py \bu) \cdot (\zb\bu)\notag\\
&\les \varp^2 \norm{\py \bu}_{\hcom})\norm{\py\bu}_{H^{m-1}_{co}}.
\end{align}
The other cases can be estimated similarly.}} Hence we can obtain the conormal estimates for $I_4,$ using the Cauchy inequality, 
\begin{align}
|I_4|\les \varp^2 \norm{\py \bu}_{\hcom}\norm{\py\bu}_{H^{m-1}_{co}}+\norm{\py\bu}_{H^{m-1}_{co}}^2
+\norm{\bu}_{\hcom}^2+\norm{\nabla p}_{H_{co}^{m-1}}\norm{\bu}_{\hcom}.\label{i4}
\end{align}

Similarly, one can obtain that 
\begin{align}
|I_5|\les \norm{\theta}_{\hcom}^2+\norm{\bu}^2_{\hcom}.\label{i5}
\end{align}
Now one can deduce \eqref{coest} from \eqref{bterm}, \eqref{i1}, \eqref{i2}, \eqref{i3}, \eqref{i6},
\eqref{i4}, \eqref{i5} and the Young inequality. The proof of Lemma \ref{le-co} is completed.
\end{proof}

{\it \textbf{Normal derivatives estimates.}}

In view of  the conormal energy estimates, we shall give the estimates of the normal derivatives $\norm{\py \bu}_{\hcoml}.$ Due to the divergence free condition, we obtain that
$$\norm{\py u_2}_{\hcoml}\les \norm{u_1}_{\hcom}\les \norm{\bu}_{\hcom}.$$
Hence it remains to estimate the normal derivatives of the tangential velocity $\norm{\py u_1}_{\hcoml}.$ 

Recall that  the vorticity $\omega$  is defined as 
$$\omega =\textup{curl}\  {\bu}=\py u_1-\px u_2.$$
For the approximate solution, we denote $\omega_a$ as
$$\omega_a=\textup{curl}\  \bu_a= \py u_{a1}-\px u_{a2}.$$
Then we can derive from  the velocity equation in \eqref{erroreq-linear}  that
\begin{align}
\pt \omega+\bu\cdot \nabla \omega_a+\bu_a\cdot \nabla \omega-\varp^2 \py^2 \omega= -\px \theta+\varp^K \textup{curl}\ R_{\bu},\label{voreq}
\end{align}
where $\textup{curl}\ {R_{\bu}}=(\py R_{u_1}-\px R_{u_2}).$ 

For the new variables 
\begin{align}\eta=\omega-\alpha u_1,\ \ \eta_a=\omega_a-\alpha u_{a1},\notag\end{align}
we have  the advantages on the boundary that
\begin{align}\eta|_{\partial \Omega}=0,\ \ \ \eta_{a}|_{\partial \Omega}=0.\label{etab}
\end{align}
Moreover, it holds that
\begin{align} \norm{\py u_1}_{\hcoml}\les \norm{\eta}_{\hcoml}+\norm{\bu}_{\hcom}.\notag\end{align}

Combining with the first equation of the velocity in \eqref{erroreq-linear}, we rewrite the vorticity equation \eqref{voreq} into the following form
\begin{align}
\pt \eta +\bu\cdot \nabla \eta_a+\bu_a\cdot\nabla \eta-\varp^2 \py^2 \eta={\alpha\px p} -\px \theta -\varp^K R_{u_1}+\varp^K \textup{curl}  R_{\bu}. \label{etaeq}
\end{align}
Taking normal derivatives on \eqref{errl-b}, one has 
\begin{align}
\pt \py \theta+\bu\cdot\nabla \py \theta_a+\bu_a\cdot\nabla \py \theta=\mathcal{C}'+\varp^K \py R_{\theta},\label{pytheta}
\end{align}
where $\mathcal{C}'=-[\py, \bu\cdot \nabla]\theta_a-[\py, \bu_a\cdot \nabla ]\theta.$

Then we have the following normal derivatives estimates.
\begin{Lemma}
For every $m\geq 1$ and every smooth solution of  the system \eqref{erroreq-linear} with partial {viscosity},  we have
\begin{align}
&\frac{d}{dt} (\norm{\eta}^2_{\hcoml}+\norm{\py \theta}^2_{\hcoml})+c_0 \varp^2 \norm{\py \eta}^2_{\hcoml} \notag\\
&\les\norm{\nabla p}_{\hcoml}\norm{\eta}_{\hcoml}+\norm{\eta}_{\hcoml}^2+\norm{\bu}_{\hcom}^2+\norm{\py\theta}_{\hcoml}^2+{{\varp^{2K-2m}}}\label{etaco}
\end{align}
for some $c_0>0.$
\end{Lemma}
\begin{proof}
We take $L^2$ inner product on equation \eqref{etaeq} with $\eta$ to get 
\begin{align}
\frac{d}{dt} \into \eta^2+\varp^2 \into |\py\eta|^2=\into \left[-\bu\cdot \nabla \eta_a \eta +{\alpha\px p} \eta +\px\theta \eta -\varp^K R_{u_1}\eta +\varp^K \textup{curl} R_{\bu}\eta\right]:=RHS,\notag
\end{align}
where $RHS$ satisfies
\begin{align}
RHS &\les \into \left[|u_1\px \eta_a \eta|+|\frac{u_2}{\varphi(y)} \varphi(y)\py \eta_a \eta|+|\px p \eta|+|\px \theta \eta|+\varp^K | R_{u_1}\eta|+{\varp^K|\curl R_\bu \eta|}\notag\right]\\
&\les \norm{\bu}^2 +\norm{\nabla p}\norm{\eta}+\norm{\eta}^2+\norm{Z (\bu, \theta)}^2+\varp^{2K-2},
\notag
\end{align}
with the aid of  the Cauchy inequality and Hardy inequality. Then we have 
\begin{align}
\frac{d}{dt} \norm{\eta}^2+\varp^2  \norm{\py\eta}^2\les \  {{  \norm{\bu}^2_{H^1_{co}}     }} +\norm{\nabla p}\norm{\eta}+\norm{\eta}^2+\norm{Z \theta}^2+\varp^{2K-2}.\label{etal2}
\end{align}
Similarly, we process $L^2$ inner product on the equation \eqref{pytheta} with $\py \theta$ to get
\begin{align}\frac{d}{dt} \norm{\py \theta}^2\les \norm{(\eta, \py\theta)}^2 + \norm{Z u}^2+\varp^{2K-2},\label{pythetal2}
\end{align}
using the Cauchy inequality and the Hardy inequality. 

We take $H_{co}^{k-1}-$inner product on the equations \eqref{etaeq} \eqref{pytheta} with $\eta$ and $\py \theta,$ respectively.  Then the desired results \eqref{etaco} can be verified by induction on $k\geq 1$. Here we omit the detail for simplicity. 
\end{proof}

{\it \textbf{Pressure estimates.}}

In this part, we  give the conormal estimates of the pressure.
In view of   the velocity equations in the partially viscous system \eqref{errl-a}, we first consider the following  system
\begin{align}
&\pt \bu-\varp^2 \py^2 \bu+\nabla p= F, \ \ \ \ \ y>0, \label{re1}\\
&\nabla\cdot \bu=0, \ \ \ \  y>0, \label{re2}\\
&u_2=0,\ \ \py u_1=\alpha u_1, \ \ \ \  y=0.\label{re3}
\end{align}
where $F$ represents the given source term. Here $F$ reads $$F=-(\bu_a\cdot \nabla \bu+\bu\cdot\nabla \bu_a)+\theta e_2+\varp^K R_\bu$$ in the linear system.

By taking divergence on the equations \eqref{re1}, using divergence free condition,  we obtain the following elliptic equation of pressure
\begin{align}
\Delta p=\nabla\cdot F, \ \ \ \ \ y>0.\label{ell}
\end{align}
One can deduce from \eqref{re2} \eqref{re3} that,  the boundary condition for \eqref{ell}, 
\begin{align}
\py p(x, 0)&=\varp^2 \py^2 u_2(x, 0)-\pt u_2(x, 0)+F_2(x, 0)\notag\\
&=-\varp^2 \py\px u_1(x, 0)+F_2(x, 0)\notag\\
&=-\varp^2\alpha \px u_1(x, 0)+F_2(x, 0).\notag
\end{align}
We can express $p$ as following 
\begin{align}p=p_1+p_2,\notag\end{align}
where $p_1$ satisfies
\begin{align}
\Delta p_1=\nabla\cdot F,  \ \ \ y>0,\ \ \ \py {p_1}(x, 0)=F_2(x, 0), \label{p1}
\end{align}
and $p_2$ solves
\begin{align}
\Delta p_2=0,\ \ \  \ y>0, \ \ \ \py p_2(x, 0)=-\alpha \varp^2 \px u_1(x, 0).\label{p2}
\end{align}
Then the desired estimates for the pressure $p$ can be achieved from the estimates of  $p_1$ and $p_2$ by standard elliptic theory.  We state the estimate of $p$ in the following 
proposition.
\begin{Proposition}
Considering the partially viscous system \eqref{re1}-\eqref{re3}, for every $m\geq 2,$ there exists $C>0$ such that for every $t\geq 0,$ we have the estimates
\begin{align}
\norm{\nabla p}_{\hcoml}\leq C(\norm{F}_{\hcoml}+\norm{\nabla\cdot F}_{H^{m-2}_{co}}+\varp^2\norm{\nabla \bu}_{\hcoml}+\norm{\bu}_{\hcoml}).\label{prop}
\end{align}
\end{Proposition}

Taking the source term $F$ in \eqref{prop} as 
\begin{align}
F=-(\bu_a\cdot \nabla \bu+\bu\cdot\nabla \bu_a)+\theta e_2+\varp^K R_\bu,\notag
 \end{align}
we can reach the desired estimates of  the pressure $p$ in the following lemma.
\begin{Lemma}\label{le-p}
For every $m\geq 2,$ and  every $\varp \in(0, 1],$ assume $(\bu, \theta)$ be  a smooth solution of \eqref{erroreq-linear} on $[0, T].$  Then it holds that
\begin{align}
\norm{\nabla p}_{\hcoml}\les (1+\varp^2)\norm{\bu}_{\hcom}+\norm{(\bu, \theta)}_{\hcoml}
+\norm{\py \bu}_{\hcoml}+\norm{\py \theta}_{H^{m-2}_{co}} +\varp^{K-m+1}.\label{pest-linear}
\end{align}
\end{Lemma}

\begin{proof}[Proof of the Theorem \ref{thmlinearco}]
Combining the {estimates} \eqref{l2est-1} \eqref{coest} \eqref{etaco} \eqref{pest-linear}, we can obtain 
\begin{align}
&\frac{d}{dt} \left(\norm{(\bu, \theta)}^2_{\hcom}+\norm{(\eta, \py \theta)}^2_{\hcoml}\right)+c_0 \varp^2 \norm{(\py \bu, \py \eta)}^2_{\hcoml} \notag\\
&\les \norm{(\bu, \theta)}_{\hcom}^2 +\norm{(\eta, \py \theta)}_{\hcoml}^2+\varp^{2K-2m}.\label{suminequ}
\end{align}
Then \eqref{thm-co} follows from  \eqref{suminequ} and the Gronwall inequality. Consequently, we have  \eqref{thm-co1} by using  \eqref{thm-co} and the Lemma \ref{le1}. The proof the Theorem \ref{thmlinearco} is completed.
\end{proof}

\section{Nonlinear Stability Estimates}
In this section, we prove the nonlinear stability  of the approximate solutions constructed in the Section $2$  under the form \eqref{app}. The zero-viscosity limit 
from the partially  viscous system \eqref{eq1} to the inviscid system \eqref{var0} will be verified with the detailed convergence rates.

In contrast  with the linear stability estimates in the Section $3,$ we shall deal with the nonlinear terms $[Z^\beta, \bold{u}\cdot \nabla]\bold{u}, [Z^\beta, \bold{u}\cdot \nabla]\theta,$ i.e.
\begin{align}
&\sum_{\gamma+\zeta=\beta, \gamma\neq 0} c_{\gamma, \zeta} Z^\gamma \bu \cdot Z^{\zeta}\nabla \bu +\bu\cdot [Z^\beta, \nabla]\bu,\ \ \ \sum_{\gamma+\zeta=\beta, \gamma\neq 0} d_{\gamma, \zeta} Z^\gamma \bu \cdot Z^{\zeta}\nabla \theta +\bu\cdot [Z^\beta, \nabla]\theta \notag
\end{align}
in the conormal energy estimates with $|\beta|\leq m$; and $[Z^\beta, \bold{u}\cdot \nabla]\eta, [Z^\beta, \bold{u}\cdot \nabla]\py\theta,$ i.e.
\begin{align}
&\sum_{\gamma+\zeta=\beta, \gamma\neq 0} \tilde{c}_{\gamma, \zeta} Z^\gamma \bu \cdot Z^{\zeta}\nabla \eta +\bu\cdot [Z^\beta, \nabla]\eta, \ \ \  \sum_{\gamma+\zeta=\beta, \gamma\neq 0} \tilde{d}_{\gamma, \zeta} Z^\gamma \bu \cdot Z^{\zeta}\nabla \py \theta +\bu\cdot [Z^\beta, \nabla]\py\theta \notag
\end{align}
in the normal derivatives estimates with $|\beta|\leq m-1$, respectively.  Here
$c_{\gamma,\zeta}, d_{\gamma,\zeta}, \tilde{c}_{\gamma,\zeta}, \tilde{d}_{\gamma, \zeta}$ are positive constants  depending on $\gamma$ and $\zeta.$
The pressure estimates should involve the nonlinear term $\bu\cdot\nabla \bu$ in the source term $F.$  In the nonlinear analysis, the estimates of $\norm{\eta}_{1, \infty}$ and $\norm{\py \theta}_{1,\infty}$ are crucial  to the closure of energy estimates. 
We deal with the term $\norm{\eta}_{1, \infty}$ in the spirit of the methods in \cite{MR} (using Maximum principle for transport-diffusion equation  and precise 
estimates for the Green's function of the operator for an approximate equation).  As for $\norm{\py \theta}_{1, \infty},$  we take advantage of a maximum principle for a particular  advective scalar  equation (cf. \cite{CC, RS}).

Let us define 
$$E_m(t)=\norm{(\bu, \theta)}_{\hcom}^2+\norm{(\eta, \py \theta)}_{\hcoml}^2+\norm{(\eta, \py \theta)}_{1,\infty}^2.$$
The proof of the Theorem \ref{maintheorem} can be reduced  to the justification of the following proposition, due to the Lemma \ref{le1}. To close the nonlinear estimates, 
the minimum expansion order is $K=6.$
\begin{Proposition}\label{priori}
Assume $(\bu, \theta)$ be a solution of \eqref{erroreq} \eqref{ini-linear} defined on $[0, T]$   with $T$  independent of $\varp$. For $m\geq 5$, $K>m,$ we have the a priori estimate
\begin{align}
E_{m}(t)\les E_m(0) +(1+t)\int_0^t (E^2_m(s)+E_m(s)) ds +\varp^{2K-2m}.
\end{align}

\end{Proposition}

Note that 
\begin{align}
\norm{\bu \cdot [Z^\beta, \nabla] \bu}\les \sum_{|\gamma|\leq m-1}\norm{u_2 \py Z^\gamma \bu}=\sum_{|\gamma|\leq m-1}\left\Vert \frac{u_2}{\varphi(y)}\varphi(y) \py Z^\gamma \bu\right\Vert\les \norm{u_2}_{W^{1, \infty}}\norm{\bu}_{\hcom}\notag
\end{align}
for $|\beta|\leq m,$ due to 
$u_2|_{y=0}=0$ and Hardy inequality.  Similarly,  we have
\begin{align}
&\norm{\bu \cdot [Z^\beta, \nabla] \theta}\les  \norm{u_2}_{W^{1, \infty}}\norm{\theta}_{\hcom}.\notag
\end{align}
By virtue of the fact that
\begin{align}\norm{Z^{\beta_1}u Z^{\beta_2} v}\les \norm{u}_{L^\infty}\norm{v}_{H^k_{co}}+\norm{v}_{L^\infty}\norm{u}_{H^k_{co}},\ \ \ |\beta_1|+|\beta_2|=k,\label{prol}\end{align}
one can deduce that
\begin{align}
&\norm{c_{\gamma, \zeta} Z^\gamma \bu \cdot Z^{\zeta}\nabla \bu}\les \norm{\nabla \bu}_{L^\infty}(\norm{\bu}_{\hcom}+\norm{\py \bu}_{\hcoml}),\\
&\norm{d_{\gamma, \zeta} Z^\gamma \bu \cdot Z^{\zeta}\nabla \theta}\les \norm{\nabla \bu}_{L^\infty}(\norm{\theta}_{\hcom}+\norm{\py \theta}_{\hcoml})+\norm{\nabla \theta}_{L^\infty}\norm{\bu}_{\hcom}.
\end{align}

Hence we can obtain the conormal energy estimates for the nonlinear system
\begin{align}
&\frac{d}{dt}(\norm{(\bu, \theta)}^2_{H_{co}^m})+c_0\varp^2\norm{\py\bu}^2_{\hcom}\notag\\
&\les (1+\norm{\bu}_{W^{1,\infty}})\left(\norm{(\theta, \bu)}^2_{\hcom}+\norm{\py( \bu, \theta)}^2_{H_{co}^{m-1}}\right)+\norm{\nabla p}_{H_{co}^{m-1}}\norm{\bu}_{\hcom}\notag\\
&\ \ \ \ +\norm{\nabla \theta}_{L^\infty}\norm{\bu}_{\hcom}^2+\varp^{2K}.\label{coest-non}
\end{align}

In the estimates of normal derivatives $\norm{(\eta, \py \theta)}_{\hcoml}$, to avoid the appearance of the terms like $\norm{\py \eta}_{L^\infty}$ and $\norm{\py \eta}_{\hcom},$  we write, for $\zeta+\gamma=\beta$, $\gamma\neq 0$, $|\beta|\leq m-1,$
\begin{align}
&\norm{\tilde{c}_{\gamma, \zeta}Z^\gamma \bu \cdot Z^{\zeta}\nabla \eta}\notag\\
&\les
\norm{Z^\gamma u_1 Z^\zeta \px \eta}+\norm{Z^\gamma u_2 Z^\zeta \py \eta}\notag\\
&\les \norm{\nabla u_1}_{L^\infty}\norm{\eta}_{\hcoml}+\norm{\eta}_{L^\infty}\norm{Z u_1}_{\hcoml}+\left\Vert\frac{1}{\varphi(y)} Z^\gamma u_2\varphi(y)Z^\zeta \py \eta\right\Vert\notag\\
&\les \norm{\bu}_{W^{1,\infty}}(\norm{\eta}_{\hcoml}+\norm{\bu}_{\hcom})
+(\norm{\bu}_{2,\infty}+\norm{Z\eta}_{L^\infty})(\norm{\eta}_{\hcoml}+\norm{\bu}_{\hcom}),\notag\\
&\les (\norm{\bu}_{W^{1,\infty}}+\norm{\bu}_{2,\infty}+\norm{Z\eta}_{L^\infty})(\norm{\eta}_{\hcoml}+\norm{\bu}_{\hcom}), \notag
\end{align}
where we have used  \eqref{prol} and divergence free condition. Indeed, for the estimate of $\norm{Z^\gamma u_2 Z^\zeta \py \eta}$, we use  the expansion of  the form
\begin{align}
\frac{1}{\varphi(y)} Z^\gamma u_2\varphi(y)Z^\zeta \py \eta=c_{\tilde{\gamma},\tilde{\zeta}}Z^{\tilde{\gamma}}\left(\frac{1}{\varphi(y)} u_2 \right)Z^{\tilde{\zeta}}(\varphi \py \eta),\label{expansion}
\end{align}
where $|\tilde{\gamma}+\tilde{\zeta}|\leq m-1$, $|\tilde{\gamma}|\neq m-1$ and $c_{\tilde{\gamma}, \tilde{\zeta}}$ is some smooth bounded coefficient (For the justification  of this expansion  \eqref{expansion},  one can refer to \cite{MR} for more details).  {{Then, with the aid of \eqref{prol} and Hardy inequality, we have 
\begin{align*}
&\left\|Z^{\tilde{\gamma}}\left(\frac{1}{\varphi(y)} u_2 \right)Z^{\tilde{\zeta}}(\varphi \py \eta)\right\|\les \left\|Z^{\tilde{\gamma}}\left(\py u_2 \right)Z^{\tilde{\zeta}}(Z \eta)\right\|\\
&\les|Z\eta|_{L^\infty}\|Z\px u_1\|_{H^{m-2}_{co}}+|Z\px u_1|_{L^\infty}\|Z\eta\|_{H^{m-2}_{co}}\\
&\les(\norm{\bu}_{2,\infty}+\norm{Z\eta}_{L^\infty})(\norm{\eta}_{\hcoml}+\norm{\bu}_{\hcom}).
\end{align*}     }}
Similarly, we can obtain that
\begin{align}
\norm{\tilde{d}_{\gamma,\zeta}Z^\gamma \bu Z^\zeta \nabla \py \theta}
\les 
& \left(\norm{\bu}_{W^{1,\infty}}+\norm{\bu}_{2, \infty}\right)\norm{\py \theta}_{\hcoml}+\norm{\py \theta}_{1, \infty}\norm{\bu}_{\hcom}.
\end{align}
Together with 
\begin{align}
\norm{\bu\cdot[Z^\beta, \nabla]\eta}\les \norm{\bu}_{W^{1,\infty}}\norm{\eta}_{\hcoml}, \ \ \ \norm{\bu\cdot[Z^\beta, \nabla]\py \theta}\les \norm{\bu}_{W^{1,\infty}}\norm{\py \theta}_{\hcoml}\notag
\end{align}
for $|\beta|\leq m-1,$  it follows that
\begin{align}
&\frac{d}{dt} (\norm{\eta}^2_{\hcoml}+\norm{\py \theta}^2_{\hcoml})+c_0 \varp^2 \norm{\py \eta}^2_{\hcoml} \notag\\
&\les\norm{\nabla p}_{\hcoml}\norm{\eta}_{\hcoml}+
 (1+\norm{\bu}_{W^{1,\infty}}+\norm{\bu}_{2,\infty}+\norm{Z\eta}_{L^\infty})(\norm{\eta}_{\hcoml}^2+\norm{\bu}_{\hcom}^2)\notag\\
 &\ \ \  +\norm{\py \theta}_{1,\infty} \norm{\bu}_{\hcom}^2+(1+\norm{\bu}_{W^{1,\infty}}+\norm{\bu}_{2,\infty})\norm{\py\theta}_{\hcoml}^2+\varp^{2K-2}.\label{etaco-non}
\end{align}

In the pressure estimates for the  nonlinear system \eqref{erroreq},  we consider the source term as 
\begin{align}
	F=-(\bu_a\cdot \nabla \bu+\bu\cdot\nabla \bu_a)+\bu\cdot\nabla \bu+\theta e_2+\varp^K R_\bu,\notag
\end{align}
 Then it follows from \eqref{prop} that
\begin{align}
\norm{\nabla p}_{\hcoml}\les \left(1+\norm{\bu}_{W^{1,\infty}}\right)\left(\norm{\bu}_{\hcom}+\norm{\py \bu}_{\hcoml}\right)+\norm{(\bu, \theta)}_{\hcoml}
+\norm{\py \theta}_{H^{m-2}_{co}} +\varp^{K-1},\label{pest-non}
\end{align}
due to
\begin{align}
&\norm{\bu\cdot \nabla \bu}_{\hcoml}\les \norm{\bu}_{W^{1,\infty}}\left(\norm{\bu}_{\hcoml}+\norm{\nabla \bu}_{\hcoml}\right),\notag\\
&\norm{\nabla \bu\cdot \nabla \bu}_{H^{m-2}_{co}}\les \norm{\nabla \bu}_{L^\infty}\norm{\nabla \bu}_{H_{co}^{m-2}}.\notag
\end{align}
Hence we can deduce from \eqref{coest-non} \eqref{etaco-non} \eqref{pest-non} that
\begin{align}
&\frac{d}{dt}\left(\norm{(\bu, \theta)}^2_{H_{co}^m}+\norm{(\eta, \py \theta)}_{\hcoml}^2\right)+c_0\varp^2\left(\norm{\py\bu}^2_{\hcom}+\norm{\py \eta}_{\hcoml}^2\right)\notag\\
&\les\left(1+\norm{\bu}_{W^{1,\infty}}+\norm{\bu}_{2,\infty}+\norm{Z\eta}_{L^\infty}\right)\left(\norm{(\theta, \bu)}^2_{\hcom}+\norm{( \eta, \py\theta)}^2_{H_{co}^{m-1}}\right)\notag\\
&\ \ \ \ +\left(\norm{\nabla \theta}_{L^\infty}+\norm{\py \theta}_{1,\infty}\right)\norm{\bu}_{\hcom}^2+\varp^{2K-2}.\label{total-1}
\end{align}

\textbf{ $L^\infty$ estimates.}

It follows from \eqref{etainfty-1} \eqref{etainfty-2} \eqref{etainfty-3} in the Lemma \ref{le2} and \eqref{coinfty} that
\begin{align}
& \norm{\bu}_{W^{1,\infty}}\les E_{m}^{1/2}(t),\ \ \norm{\bu}_{2,\infty}\les E_{m}^{1/2}(t),\ \ \norm{\nabla \bu}_{1, \infty}\les  E_{m}^{1/2}(t), \ \ m\geq m_0+3,\label{inf-u}\\
&\norm{\theta}_{1,\infty}\les \norm{\theta}_{m}+\norm{\py \theta}_{m-1}\les E_m^{1/2}(t), \ \ \norm{\theta}_{2,\infty}\les E_m^{1/2}(t), \ \ m\geq m_0+3, \label{inf-theta}
\end{align}
for $m_0\geq 1.$

Combining   the inequalities  \eqref{total-1}  \eqref{inf-u} and \eqref{inf-theta}, we still need to estimate
$$\norm{\eta}_{1, \infty},\ \ \norm{\py \theta}_{1, \infty}$$
to close the energy estimates. 
 Similar to the  Proposition 13 in \cite{MR} on  the estimates of $\norm{\eta}_{1,\infty},$ we can derive the following results:
\begin{Lemma}
For $m\geq 5,$ we have the estimate
\begin{align}\norm{\eta}^2_{1,\infty}\les E_m(0)+(1+t)\int_0^t (E_m^2(s)+E_m(s))ds+\varp^{2K-2m}.\label{eta-1inf}\end{align}
\end{Lemma}
\begin{proof}
{{Recall that $\eta$ satisfies the following nonlinear equations
\begin{align*}
&\pt \eta +(\bu_a+\bu)\cdot\nabla \eta-\varp^2 \py^2 \eta=\mathscr{R}:=-\bu\cdot \nabla \eta_a+\alpha\px p-\px \theta -\varp^K R_{u_1}+\varp^K \textup{curl}  R_{\bu},\\
&\eta=0\ \ \ {\rm on}\ \ \ \{y=0\}.\end{align*}
Here $\eta$ solves a transport-diffusion equation and satisfies a homogenous Dirichlet boundary condition.  The estimates of $\norm{\eta}_{L^\infty}$ and $\norm{Z_1\eta}_{L^\infty}$
can be derived by maximal principle for transport-diffusion equation. The estimate of $\norm{Z_2 \eta}_{L^\infty}$ would be the  most difficult one, due to the commutators between $Z_3$ and $\py^2$.  Inspired by the idea in \cite{MR}, one can similarly use the estimates on the Green's function of an approximate equation for  the above equation near the boundary, i.e.
\begin{align*}
&\pt \eta +(u_{1a}(t, x, 0)+u_1(t, x, 0))\px \eta+y\py(u_{2a}(t, x, 0)+u_2(t, x, 0))\py \eta-\varp^2 \py^2 \eta\\
&=\mathscr{R}-G,\tag{*}
\end{align*}
where $G=[u_{1a}+u_1-u_{1a}(1, x, 0)-u_{1}(t, x, 0)]\px \eta+[u_{2a}+u_2-y\py(u_{2a}(t, x, 0)-u_2(t, x, 0)]\py\eta.
$
Notice that the solution of the equation $(*)$ can be expressed by the generator $S(t, \tau)$ of the operator in the left hand side of $(*)$ through Duhamel formula. 
Namely, for $\forall t\geq \tau,$
\begin{align*}
\eta(t)=S(t, \tau)\eta_0+\int_0^t S(t, \tau)(\mathscr{R}-G)(\tau)d\tau.
\end{align*}
This,  together  with 
\begin{align*}
\norm{y\py S(t, \tau) \eta_0}_{L^\infty}\les \norm{\eta_0}_{L^\infty}+\norm{y\py \eta_0}_{L^\infty},
\end{align*}
gives  the estimate of $\norm{Z_2\eta}_{L^\infty},$
$$\norm{Z_2 \eta}_{L^\infty}\leq\left(\norm{\eta_0}_{1, \infty}+\int_0^t \norm{\mathscr{R}-G}_{1,\infty}\right).$$ Then \eqref{eta-1inf} can be deduced from  a closed conormal argument for $\norm{(\mathscr{R}-G)}_{1,\infty}\les E_m^2(t)+E_m(t), m\geq 5$.  Here we omit the details for simplicity. }}

\end{proof}

Therefore, it remains to estimate $\norm{\py \theta}_{1,\infty}$, which can be derived from  the  maximum principle for an advective scalar  equation (cf. \cite{CC, RS}).
\begin{Lemma}
For $m\geq 5,$ we have the estimate
\begin{align}\norm{\py \theta}^2_{1,\infty}\les E_m(0)+ \int_0^t E_m(s)  ds+\varp^{2K-4}.\label{pytheta-1inf}\end{align}
\end{Lemma}
\begin{proof}
We rewrite \eqref{err-b} in the following form
\begin{align}
\pt \theta +(\bu+\bu_a)\cdot \nabla \theta=f(t,x,y):=-\bu\cdot \nabla \theta_a+\varp^K R_\theta,\ \ y>0,
\end{align}
where $\bu, \bu_a$ satisfy 
$$\nabla\cdot \bu=0,\ \ \nabla\cdot \bu_a=0, \ \ u_2(t,x,0)=0,\ \ u_{a2}(t, x, 0)=0.$$
We transform the problem into the whole space by defining $\tilde{\theta}, \tilde{\theta}_a, \tilde{\bu}, \tilde{\bu}_a, \tilde{f}:$
\begin{align}
&(\tilde{\theta},\tilde{\theta}_a, \tilde{\bu},\tilde{\bu}_a)=(\theta, \theta_a,  \bu, \bu_a)(t, x, y), \ y>0;\notag\\
&{{(\tilde{\theta}, \tilde{\theta}_a, \tilde{\bu}, \tilde{\bu}_a)=-(\theta, \theta_a,  (-u_1, u_2), (-u_{1a}, u_{2a}))(t, x, -y),\  y<0}};\notag\\
&\tilde{f}=-\tilde{\bu}\cdot\nabla \tilde{\theta}_a+\varp^K \tilde{R}_{{\theta}}.\notag
\end{align}
Then one has 
\begin{align}
\pt \tilde{\theta} +(\tilde{\bu}+\tilde{\bu}_a)\cdot \nabla \tilde{\theta}=\tilde{f},\ \ \ \tilde{\theta}_0= \tilde{\theta}|_{t=0}
\end{align}
with $\nabla\cdot \tilde{\bu}=0$ and  $\nabla\cdot \tilde{\bu}_a=0.$
Applying maximum principle for advective  scalar equations (see \cite{RS}, $3.2$ A priori Bounds. and \cite{CP, CC}), we obtain that, for $t\in [0, T],$
\begin{align}
\norm{\tilde{\theta}}_{L^\infty}&\les \norm{\tilde{\theta}|_{t=0}}_{L^\infty}+\int_0^t \norm{\tilde{f}(\tau)}_{L^\infty}d\tau\notag\\
&\les \varp^{K+1} \norm{\theta_0}_{L^\infty}+\int_0^t \norm{\bu(\tau)}_{L^\infty}d\tau+\varp^{K}.
\end{align}
Similarly, we can derive the $L^\infty$-estimates of $\nabla\tilde{\theta}$ and $Z\py\tilde{\theta},$
\begin{align}
\norm{\nabla \tilde{\theta}}_{L^\infty}&\les  \varp^{K+1} \norm{\nabla\theta_0}_{L^\infty}+\int_0^t  (\norm{\bu(\tau)}_{L^\infty}+\norm{\nabla \bu(\tau)}_{L^\infty}) d\tau+\varp^{K-1},\\
\norm{Z\py \tilde{\theta}}_{L^\infty}&\les  \varp^{K+1} \norm{Z \py\theta_0}_{L^\infty}+\int_0^t  (\norm{\bu(\tau)}_{W^{1,\infty}}+\norm{\bu}_{1,\infty}+\norm{Z\py \bu(\tau)}_{L^\infty}) d\tau+\varp^{K-2},
\end{align}
due to 
\begin{align}
\norm{ u_2\py \py\theta_a}_{L^{\infty}}\leq \left\Vert\frac{u_2}{\varphi(y)}\varphi(y)\py \py \theta_a\right\Vert_{L^{\infty}}\les \norm{u_2}_{W^{1,\infty}}.
\end{align}
Hence the proof of \eqref{pytheta-1inf} is completed with the aid of Lemma \ref{le2}.
\end{proof}
Therefore, the proof of Proposition \ref{priori} can be obtained by combining \eqref{total-1} \eqref{eta-1inf} \eqref{pytheta-1inf} and  the standard continuous argument.

\textbf{Acknowledgements.}
{The research of the authors  was partially  supported by the National Natural Science Foundation of China (NSFC) under grants 12101350, 12171267 and 12271284;  China Postdoctoral  Science Foundation under grant 2021M691818;  Natural  Science Foundation of Jiangsu Province under grant  BK20220792. The authors are  grateful to the anonymous referees for their  careful reading of the manuscript and the numerous
very helpful suggestions which have helped to clarify some important points and improve
the exposition of this  paper greatly .}

\textbf{Data availability.}
Data sharing not applicable to this article as no datasets were generated or analysed during the current study.

\textbf{Conflict of interest.}
The authors declare that they have no conflict of interest.

\end{document}